\documentclass[11pt]{amsart}

\usepackage[T1]{fontenc}
\usepackage{amsmath,amssymb,amsthm}
\usepackage{graphicx}
\usepackage[margin=1in]{geometry}
\usepackage[hidelinks]{hyperref}
\hypersetup{
  pdftitle={Polynomial bound for the localization length of the Lorentz mirror model on the 1D cylinder},
  pdfauthor={Linjun Li},
  pdfkeywords={Lorentz lattice gas, Manhattan pinball, localization, cylinder}
}

\title[Lorentz mirror on a cylinder]
{Polynomial bound for the localization length of the Lorentz mirror model on the 1D cylinder}
\author[L. Li]{Linjun Li}
\address{Department of Mathematics, University of Pennsylvania,
Philadelphia, PA}
\email{linjun@sas.upenn.edu}
\date{}

\newtheorem{theorem}{Theorem}[section]
\newtheorem{lemma}[theorem]{Lemma}
\newtheorem{prop}[theorem]{Proposition}
\newtheorem{cor}[theorem]{Corollary}
\theoremstyle{definition}

\newcommand{\Z}{\mathbb Z}
\newcommand{\Prob}{\mathbb P}
\newcommand{\E}{\mathbb E}
\newcommand{\cC}{\mathcal C}
\newcommand{\Cross}{\mathrm{Cross}}
\newcommand{\Col}{\mathrm{Col}}
\newcommand{\Wind}{\mathrm{Wind}}
\newcommand{\Trav}{\mathrm{Trav}}
\newcommand{\Block}{\mathrm{Block}}
\newcommand{\Plane}{\mathrm{PlaneLoc}}

\begin{document}

\begin{abstract}
We prove polynomial upper bounds for the localization length of the Lorentz mirror model and the Manhattan model on the even cylinder.  We first show that a fixed positive lower bound for short-direction crossings of a $100n\times n$ rectangle implies localization on scale $O(n^{10})$.  The proof is genuinely cylindrical and combines winding barriers with a two-site switching and double-counting argument.  Together with a planar confinement argument proved here, this yields unconditional cylinder localization for both models.  For Lorentz mirrors, a planar escape estimate ensures the required crossing lower bound; for Manhattan mirrors, planar confinement handles any fixed scale at which the crossing lower bound fails.
\end{abstract}

\maketitle

\section{Introduction}\label{sec:intro}

Random mirror models are deterministic dynamics in a random environment: a ray
on the square lattice is transmitted or reflected by a diagonal mirror at each
vertex.  We prove polynomial localization bounds on the even cylinder for both
the Lorentz mirror and Manhattan models.

Our main contribution is a genuinely cylindrical argument.  A lower bound for
short-direction crossings of planar rectangles produces many closed trajectories
winding around the cylinder.  These trajectories are topological barriers, and a
two-site switching and double-counting argument then suppresses trajectories that
cross a long cylindrical block.  The unconditional consequences differ between
the models: planar escape estimates force the required crossing bound for Lorentz
mirrors, whereas for Manhattan mirrors a low crossing probability at one fixed
scale may instead yield planar confinement.

Let
\[
        \cC_n=\Z\times(\Z/2n\Z)
\]
be the cylinder of circumference $2n$, with second coordinate modulo $2n$.
Write $\Prob_p^{\mathrm L}$ and $\Prob_p^{\mathrm M}$ for the Lorentz and
Manhattan product laws on $\Z^2$ or $\cC_n$, as context indicates, and $\Prob_p$
when either fixed model is meant.

\emph{Lorentz mirror model.}
Fix $p\in(0,1)$.  Independently at each vertex, place no mirror with probability
$1-p$, and a north-west (NW) or north-east (NE) mirror with probability $p/2$
each.

\begin{samepage}
\emph{Manhattan model.}
On $\Z^2$, orient horizontal edges eastward on odd rows and westward on even
rows, and orient vertical edges northward on even columns and southward on odd
columns.  Since the circumference is $2n$, these directions descend to $\cC_n$.
Independently place a mirror with probability $p$ and no mirror otherwise; at
$u=(x,y)$ it is NW when $x-y$ is even and NE otherwise.  The same local rules
define both models on $\Z^2$ and $\cC_n$.
\end{samepage}

\begin{figure}
    \centering
    \includegraphics[width=0.6\linewidth]{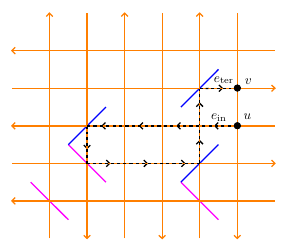}
    \caption{The pink and blue segments illustrate NW and NE mirrors,
    respectively.  The yellow arrows illustrate the underlying Manhattan
    directions, and the dashed line illustrates a trajectory.  The labels
    $e_{\rm in}$ and $e_{\rm ter}$ mark its initial and terminal directed
    edges, and $u,v$ mark the indicated endpoint vertices.}
    \label{fig:model}
\end{figure}

A ray follows lattice edges, reflecting at mirrors and otherwise going straight;
see Figure~\ref{fig:model}.  A \emph{directed edge} is an ordered nearest-neighbor
edge $(u,v)$ and, with the environment, determines the trajectory.  Local rules
are reversible.  We sometimes count only Manhattan-oriented trajectories, since
the reverse orientation represents the same geometric trajectory.

For $M\ge1$, define the cylinder-localization event
\[
        B_n^{(M)}
        =
        \left\{
        \begin{array}{l}
        \text{every trajectory whose initial vertex lies in }
        \{0\}\times(\Z/2n\Z)\\
        \text{remains in }[-M,M]\times(\Z/2n\Z)
        \end{array}
        \right\}.
\]
We also use a planar rectangle-crossing event.  For $m,h\ge1$, let
\[
        R_{m,h}=\{1,\ldots,m\}\times\{1,\ldots,h\}.
\]
For $1\le i,j\le m$, let $\Cross_{i,j}(m,h)$ be the event that there exists a
trajectory segment which enters $R_{m,h}$ through $((i,0),(i,1))$, exits through
$((j,h),(j,h+1))$, and has all intermediate vertices in $R_{m,h}$.  Set
\[
        \Cross(m,h)=\bigcup_{1\le i,j\le m}\Cross_{i,j}(m,h),
        \qquad
        \Col(m,h)=\bigcup_{1\le i\le m}\Cross_{i,i}(m,h).
\]
Thus $\Cross$ is a vertical crossing event and $\Col$ also requires equal
entrance and exit columns.  Translations include the boundary edges; all are
allowed for Lorentz mirrors, while Manhattan translations preserve parity.

Throughout the main statements, $c_0>0$ denotes a universal constant, chosen
small enough that Proposition~\ref{prop:planar} holds.

\paragraph{Lorentz mirror model.}

\begin{theorem}[Cylinder localization from crossings: Lorentz model]
\label{thm:lorentz-one-scale}
Fix $p\in(0,1)$.  There are constants $C=C(p)<\infty$ and
$n_*=n_*(p)<\infty$ such that, whenever $n\ge n_*$ satisfies
\[
        \Prob_p^{\mathrm L}\bigl(\Cross(100n,n)\bigr)\ge c_0,
\]
then, for every $M\ge Cn^{10}$,
\[
        \Prob_p^{\mathrm L}\bigl(B_n^{(M)}\bigr)
        \ge
        1-C\exp\left\{-\frac{M}{Cn^{10}}\right\}.
\]
\end{theorem}

\paragraph{Manhattan model.}

The same crossing assumption gives a cylinder estimate for the Manhattan model,
with a different local switching operation in the proof.

\begin{theorem}[Cylinder localization from crossings: Manhattan model]
\label{thm:manhattan-one-scale}
Fix $p\in(0,1)$.  There are constants $C=C(p)<\infty$ and
$n_*=n_*(p)<\infty$ such that, whenever $n\ge n_*$ and
$\Prob_p^{\mathrm M}(\Cross(100n,n))\ge c_0$, then, for every $M\ge Cn^{10}$,
\[
        \Prob_p^{\mathrm M}\bigl(B_n^{(M)}\bigr)
        \ge
        1-C\exp\left\{-\frac{M}{Cn^{10}}\right\}.
\]
\end{theorem}

\emph{Even circumference.}
The even circumference is essential.  It is needed for the Manhattan
orientation to descend to the cylinder, while odd circumferences in the
Lorentz model admit infinite trajectories.

Theorems~\ref{thm:lorentz-one-scale} and
\ref{thm:manhattan-one-scale} have a common geometric proof, except for the
model-specific local switching step.  The following result is instead planar
and does not use the topology of the cylinder.

\begin{prop}[Low crossing implies planar confinement]\label{prop:planar}
Fix $p\in(0,1)$ and either one of the two models.  Suppose that for some
$n_0>100$,
\[
        \Prob_p\bigl(\Cross(100n_0,n_0)\bigr)\le c_0.
\]
Let $\Plane_m$ be the event that every planar trajectory whose initial vertex
lies in $[-1,1]^2$ remains in $[-m,m]^2$.  There are universal constants
$c,C>0$ such that
\[
        \Prob_p(\Plane_m)
        \ge 1-\exp\{-cm/n_0\},
        \qquad m\ge Cn_0.
\]
\end{prop}

We now record the unconditional consequences.  For Lorentz mirrors, let $E_m$
be the event that a fixed ray from the origin reaches outside $[-m,m]^2$.
Kozma and Sidoravicius proved that
$\Prob_p^{\mathrm L}(E_m)\ge(2m+1)^{-1}$~\cite{kozma2015lower}.  The hypothesis
of Proposition~\ref{prop:planar} cannot therefore hold: at a scale $n_0$ it
would imply $E_m\subseteq\Plane_m^c$ and
\[
        \frac1{2m+1}
        \le \Prob_p^{\mathrm L}(E_m)
        \le \Prob_p^{\mathrm L}(\Plane_m^c)
        \le e^{-cm/n_0},
\]
a contradiction for all sufficiently large $m$.  Hence
\[
        \Prob_p^{\mathrm L}\bigl(\Cross(100n,n)\bigr)>c_0
        \qquad\text{for every }n>100,
\]
and Theorem~\ref{thm:lorentz-one-scale} applies whenever in addition
$n\ge n_*(p)$.

\begin{cor}[Unconditional Lorentz cylinder localization]\label{cor:lorentz}
For the Lorentz mirror model and every fixed $p\in(0,1)$, there are constants
$C=C(p)<\infty$ and $n_*=n_*(p)<\infty$ such that, for all $n\ge n_*$ and
$M\ge Cn^{10}$,
\[
        \Prob_p^{\mathrm L}\bigl(B_n^{(M)}\bigr)
        \ge
        1-C\exp\left\{-\frac{M}{Cn^{10}}\right\}.
\]
\end{cor}

For Manhattan mirrors, Theorem~\ref{thm:manhattan-one-scale} applies at every
sufficiently large width if its crossing inequality holds at all such widths.
Otherwise Proposition~\ref{prop:planar} at one fixed low-crossing scale,
followed by a plane-to-cylinder comparison, gives the same bound at all large
target widths.

\begin{cor}[Unconditional Manhattan cylinder localization]\label{cor:manhattan}
For the Manhattan model and every fixed $p\in(0,1)$, there are constants
$C=C(p)<\infty$ and $n_*=n_*(p)<\infty$ such that, for all $n\ge n_*$ and
$M\ge Cn^{10}$,
\[
        \Prob_p^{\mathrm M}\bigl(B_n^{(M)}\bigr)
        \ge
        1-C\exp\left\{-\frac{M}{Cn^{10}}\right\}.
\]
\end{cor}

We do not expect the exponent $10$ to be optimal; see
\cite[Page~238]{spencer2012duality}.  Under the crossing assumptions above, the
natural conjectural exponent is $2$.  When $p\le Cn^{-1}$, the localization
length is known to be $O(p^{-2})$~\cite{ryan2020manhattan}.  In a different
direction, for $d\ge4$ and sufficiently small $p$, Lorentz mirror trajectories
remain diffusive up to every fixed polynomial time scale in
$p^{-1}$~\cite{elboim2025diffusivity}.  The Manhattan model is also a network
model for a unitary quantum dynamics; see
\cite{li2021manhattan,spencer2012duality,cardy2010quantum,beamond2003quantum}.
For background on Anderson localization and related localization and
delocalization results for random operators and band matrices, see
\cite{anderson1958absence,klein1990localization,bourgain2013lower,
binder2015fluctuations,Li2D2022,LiThesis2022,LiZhang3D,LiSupport2026,
dubova2025delocalization,yau2025delocalization,drogin2025localization}.

Here is the proof architecture.  Under the assumed lower bound for
$\Cross(100n,n)$, Lemma~\ref{lem:glue} produces a same-column crossing of a
$100n\times2n$ rectangle with probability at least $c(p,c_0)n^{-4}$.
Lemma~\ref{lem:column-to-wind} projects this crossing to a winding closed
trajectory, and Lemma~\ref{lem:sep} shows that its vertex set is a separating
barrier.
In a length-$N\asymp n^{10}$ strip, on the order of $n^9$ independent trials give at least
$c(p,c_0)n^5$ barriers with high probability.  Every horizontal traversal meets
them; a recoverable two-site switch and double counting make one-block traversal
unlikely.  Independent blocks yield the exponential tail.  The planar
proposition follows separately by coarse graining.

\section{Crossings and gluing}\label{sec:cross}

This section proves the elementary gluing estimate used to create winding
trajectories.  The estimate is stated in planar rectangles; later, in
Section~\ref{sec:winding}, a same-column crossing of height $2n$ will be
projected to the cylinder and become a winding closed trajectory.

We shall use the following language throughout the paper.  A trajectory segment
is a finite sequence of directed edges
\[
        e_k=(v_{k-1},v_k)
\]
which follows the local scattering rule at every intermediate vertex.  If
$A=(e_1,\ldots,e_r)$ is such a directed-edge sequence, define its inverse by
\[
        A^{-1}=(\bar e_r,\ldots,\bar e_1),
        \qquad \bar e=(v,u)\ \hbox{for }e=(u,v).
\]
If the terminal vertex of $A$ is the initial vertex of another directed-edge
sequence $B$, then $AB$ denotes their concatenation.  This concatenation is a
trajectory segment whenever the local scattering rule also holds at the joining
vertex.  A
\emph{marked occurrence} of a path at a vertex means the incoming half-edge and
the outgoing half-edge used during that visit.  This distinction matters when a
trajectory visits the same vertex more than once.  In a finite set, an exiting
trajectory segment cannot repeat a directed edge before exiting: a first repeat,
together with reversibility, would force the preceding directed edge also to
repeat.  Consequently an exiting segment has at most four marked occurrences
at any fixed vertex, one for each possible incoming directed edge.

\begin{lemma}[Gluing Lemma]\label{lem:glue}
For either model and all $m\ge1$ and $n\ge2$,
\[
        \Prob_p\bigl(\Col(m,2n)\bigr)
        \ge
        (1-p)^2m^{-4}\Prob_p\bigl(\Cross(m,n)\bigr)^2 .
\]
\end{lemma}

\begin{proof}
Let
\[
        q=\Prob_p\bigl(\Cross(m,n)\bigr).
\]
Since $\Cross(m,n)$ is the union of $m^2$ events $\Cross_{i,j}(m,n)$, there are
columns $i,j$ such that
\[
        \Prob_p\bigl(\Cross_{i,j}(m,n)\bigr)\ge q\,m^{-2}.
\]

We first treat the Lorentz model, and also the Manhattan model when $n$ is even.
Let $\rho(x,y)=(x,n+1-y)$ be reflection in the horizontal line
$y=(n+1)/2$.  For Lorentz mirrors, $\rho$ interchanges the equally likely NW
and NE pairings.  For Manhattan mirrors,
$x-(n+1-y)\equiv x-y+1\pmod2$ when $n$ is even, so $\rho$ interchanges the
parity classes exactly as it interchanges the mirror pairings.  Thus in both
models reflection preserves the law, and reversal of the reflected segment
maps an $i$-to-$j$ crossing to a $j$-to-$i$ crossing.  Hence
\[
        \Prob_p\bigl(\Cross_{j,i}(m,n)\bigr)
        =
        \Prob_p\bigl(\Cross_{i,j}(m,n)\bigr).
\]
The vertical shift by $n$ is parity preserving.  Hence an $i$-to-$j$ lower
crossing and an independent $j$-to-$i$ upper crossing, using disjoint site
variables, occur with probability at least $q^2m^{-4}$.  The two segments share
the directed interface edge $((j,n),(j,n+1))$.  Retain this edge as the last edge
of the lower segment, delete its duplicate as the first edge of the upper
segment, and concatenate the remaining directed-edge sequences.  This gives a
same-column crossing of height $2n$.  Hence in these cases
\[
        \Prob_p\bigl(\Col(m,2n)\bigr)\ge q^2m^{-4},
\]
which is stronger than the claimed bound.

It remains only to handle the Manhattan model when $n$ is odd.  Any crossing of
$R_{m,n}$ contains, by stopping at the first edge crossing from height $n-1$ to
height $n$, an initial crossing of $R_{m,n-1}$.  Therefore
\[
        \Prob_p\bigl(\Cross(m,n-1)\bigr)\ge q .
\]
Since $n-1$ is even, the preceding argument applied at height $n-1$ gives
\[
        \Prob_p\bigl(\Col(m,2n-2)\bigr)\ge q^2m^{-4}.
\]
On the event $\Col(m,2n-2)$ choose, by a fixed measurable rule, the least column
$i$ witnessing a crossing.  The event and witness use only site variables up
to height $2n-2$; those at $(i,2n-1)$ and $(i,2n)$ remain independent and
unrevealed.  With conditional probability $(1-p)^2$ both are vacant, extending
the ray straight to height $2n$.  Thus
\[
        \Prob_p\bigl(\Col(m,2n)\bigr)
        \ge
        (1-p)^2 q^2 m^{-4},
\]
as required.
\end{proof}

\section{Winding barriers}\label{sec:winding}

For integers $a<b$ write
\[
        T^{(n)}_{a,b}=\{a+1,a+2,\ldots,b\}\times(\Z/2n\Z).
\]
A closed nearest-neighbor path $W=(w_0,w_1,\ldots,w_r)$ on $\cC_n$, with
$w_r=w_0$, has \emph{vertical winding number} $\mathrm{wind}(W)=k$ if, after
choosing a lift $\widetilde W=(\widetilde w_0,\ldots,\widetilde w_r)$ to
$\Z^2$, one has
\[
        \widetilde w_r=\widetilde w_0+(0,2nk).
\]
The integer $k$ is independent of the chosen lift and changes sign under reversal
of the orientation of $W$.  We call $W$ \emph{winding} if $k\ne0$.  If a closed
trajectory visits a vertex more than once, then each visit is understood as a
separate local passage through that vertex.

We shall use the following elementary topological fact in a purely discrete form.
A set $S\subset T^{(n)}_{a,b}$ is said to separate the two sides of the strip if
every nearest-neighbor path in $\cC_n$ from $\{x\le a\}$ to $\{x\ge b+1\}$ meets
$S$.

\begin{lemma}\label{lem:sep}
Let $W$ be a closed nearest-neighbor path on $\cC_n$ with non-zero vertical
winding number and with $V(W)\subset T^{(n)}_{a,b}$.  Then $V(W)$ separates the
two sides of the strip $T^{(n)}_{a,b}$.
\end{lemma}

\begin{proof}
It suffices to prove the statement for a simple closed subpath of non-zero
winding number.  Indeed, decompose the closed walk $W$ by successively cutting at
the first repeated vertex.  This writes $W$ as a concatenation of simple closed
cycles, and the vertical winding number is the sum of the winding numbers of
these cycles.  Since $\mathrm{wind}(W)\ne0$, at least one cycle, say $\Gamma$, has
non-zero winding number.  Also $V(\Gamma)\subset V(W)$.

Realize $\cC_n$ as the annulus $\mathbb R\times(\mathbb R/2n\Z)$ and realize
lattice paths as polygonal curves.  Since $\Gamma$ is simple and has non-zero
winding number, it is an essential simple closed curve in the finite annulus
\[
        [a+1/2,b+1/2]\times(\mathbb R/2n\Z).
\]
By the Jordan curve theorem on the annulus, an essential simple closed curve
separates the two boundary components of this annulus.  Therefore every
continuous curve from the left boundary $\{a+1/2\}\times(\mathbb R/2n\Z)$ to the
right boundary $\{b+1/2\}\times(\mathbb R/2n\Z)$ intersects $\Gamma$.

Now let $\gamma$ be a nearest-neighbor path in $\cC_n$ from $\{x\le a\}$ to
$\{x\ge b+1\}$, and consider its polygonal realization.  Let $t_R$ be its first
crossing of $\{b+1/2\}\times(\mathbb R/2n\Z)$, and let $t_L<t_R$ be its last
crossing of $\{a+1/2\}\times(\mathbb R/2n\Z)$ before time $t_R$.  The subcurve
$\gamma|_{[t_L,t_R]}$ lies in the annulus and joins its boundary components,
so it intersects $\Gamma$.  Since both curves use lattice edges, the
intersection contains a common lattice vertex.  Thus $\gamma$ meets
$V(\Gamma)\subset V(W)$.
\end{proof}

We next formalize the passage from same-column crossings to winding barriers.
For $\ell\ge1$, let
\[
        I_\ell=\{200\ell n+1,\ldots,200\ell n+100n\},
        \qquad S_\ell=T^{(n)}_{200\ell n,\,200\ell n+100n}.
\]
The map
\[
 \pi_\ell:R_{100n,2n}\longrightarrow S_\ell,
 \qquad
 \pi_\ell(x,y)=(x+200\ell n,[y]_{2n}),
\]
is a bijection of sites.  It preserves the Manhattan parity of $x-y$ because
both $200\ell n$ and $2n$ are even.  We transport the event
$\Col(100n,2n)$ to the cylinder through $\pi_\ell$ and denote the resulting
event by $\Wind_\ell$.  Choose a witnessing crossing deterministically.  Its
first and last directed edges have the same cylinder projection.  Delete the
last, repeated edge and call the resulting closed trajectory $W_\ell$.  In the
Manhattan model give it its Manhattan orientation, fixed from now on.

\begin{lemma}\label{lem:column-to-wind}
On $\Wind_\ell$, $W_\ell$ is a closed trajectory on $\cC_n$.  It satisfies
$|\mathrm{wind}(W_\ell)|=1$ and $V(W_\ell)\subset S_\ell$.  Consequently
$V(W_\ell)$ separates $\{x\le200\ell n\}$ from
$\{x\ge200\ell n+100n+1\}$.
\end{lemma}

\begin{proof}
The initial edge $((i,0),(i,1))$ and terminal edge
$((i,2n),(i,2n+1))$ have the same projection.  After deleting the latter, a
lift of the resulting closed path runs from $(i,0)$ to $(i,2n)$, so its winding
number is $1$ in one orientation and $-1$ in the other.  Its first coordinate
lies in $I_\ell$, hence $V(W_\ell)\subset S_\ell$; now apply
Lemma~\ref{lem:sep}.
\end{proof}

\begin{lemma}\label{lem:many}
Assume $\Prob_p(\Cross(100n,n))\ge c_0$.  Let $C_1\ge1$ and
$N=\lfloor C_1n^{10}\rfloor$ and
$L=\lfloor C_1n^9/1000\rfloor$.  For all sufficiently large $n$,
$S_1,\ldots,S_L\subset T^{(n)}_{0,N}$.  Moreover, the events
$\Wind_1,\ldots,\Wind_L$ are independent, and there are constants
$c_1,C_2>0$, depending only on $p$ and $c_0$, such that
\[
        \Prob_p\left(\#\{1\le\ell\le L:\Wind_\ell\hbox{ occurs}\}
        \ge c_1C_1n^5\right)
        \ge 1-\frac{C_2}{C_1n} .
\]
\end{lemma}

\begin{proof}
For large $n$, $200Ln+100n\le(C_1/5)n^{10}+100n<N$, proving the slice
containment.  The event $\Wind_\ell$ is
determined only by the mirror variables in the slice
$S_\ell$, and the slices $S_1,\ldots,S_L$ are disjoint.  Hence these events are
independent.  By Lemma~\ref{lem:glue}, for each $\ell$,
\[
        \Prob_p(\Wind_\ell)
        \ge (1-p)^2(100n)^{-4}\Prob_p(\Cross(100n,n))^2
        \ge c_w n^{-4},
        \qquad c_w=100^{-4}(1-p)^2c_0^2 .
\]
Let $X=\sum_{\ell=1}^L {\bf 1}_{\Wind_\ell}$.  For all sufficiently large $n$,
\[
        \E X\ge \frac{c_wC_1}{2000}n^5 .
\]
Since $X$ is a sum of independent Bernoulli random variables,
$\operatorname{Var}X\le \E X$.  Chebyshev's inequality gives
\[
        \Prob_p\left(X<\frac12\E X\right)
        \le \frac{4}{\E X}
        \le \frac{C_2}{C_1n},
\]
a weaker form than the natural $O((C_1n^5)^{-1})$ bound but sufficient here.
Taking $c_1=c_w/4000$ proves the lemma.
\end{proof}

\section{Surgery and the one-block estimate}\label{sec:surgery}

Fix $N$ and work in the finite strip $D_N=T^{(n)}_{0,N}$.  We use the following
boundary directed edges:
\[
\begin{aligned}
\mathcal E_l^+&=\{((0,y),(1,y)):y\in\Z/2n\Z\}, &
\mathcal E_l^-&=\{((1,y),(0,y)):y\in\Z/2n\Z\},\\
\mathcal E_r^+&=\{((N,y),(N+1,y)):y\in\Z/2n\Z\}, &
\mathcal E_r^-&=\{((N+1,y),(N,y)):y\in\Z/2n\Z\}.
\end{aligned}
\]
In the Manhattan model, whenever we count Manhattan-oriented trajectories, these
sets are restricted to the directed edges compatible with the Manhattan
orientation: with the convention fixed above, $y$ is odd for
$\mathcal E_l^+$ and $\mathcal E_r^+$, and $y$ is even for
$\mathcal E_l^-$ and $\mathcal E_r^-$.  This restriction loses no geometric
traversal, since the inverse of every non-Manhattan-oriented traversal is
Manhattan-oriented.

A trajectory segment with all interior vertices in $D_N$ is called left--right,
right--left, left--left, or right--right according as its initial and terminal
edges lie in
\[
\mathcal E_l^+\times\mathcal E_r^+,
\quad
\mathcal E_r^-\times\mathcal E_l^-,
\quad
\mathcal E_l^+\times\mathcal E_l^-,
\quad
\mathcal E_r^-\times\mathcal E_r^+,
\]
respectively.  Let $\Trav_N$ be the event that a left--right or right--left
traversal exists.

We first record two deterministic facts about trajectories in a finite strip.

\begin{lemma}[Exit map]\label{lem:exit}
For every configuration in $D_N$, each inward boundary edge in
$\mathcal E_l^+\cup\mathcal E_r^-$ determines a unique trajectory segment whose
terminal edge lies in $\mathcal E_l^-\cup\mathcal E_r^+$.  The resulting exit map
\[
        F:\mathcal E_l^+\cup\mathcal E_r^-
          \longrightarrow \mathcal E_l^-\cup\mathcal E_r^+
\]
is a bijection.  The same statement holds in the Manhattan model after the
Manhattan-orientation restriction described above.
\end{lemma}

\begin{proof}
Starting from an inward boundary edge, follow the deterministic local rule.  If
the trajectory did not exit $D_N$, then after its first interior edge it would
move forever in the finite set of directed edges with both endpoints in $D_N$;
hence it would eventually enter a directed cycle.  Let $e_t$ be the first edge of
that cycle reached by the trajectory.  The predecessor of $e_t$ inside the cycle
and the predecessor $e_{t-1}$ by which the trajectory first reached the cycle are
two predecessors of the same directed edge.  This is impossible because the local
scattering map is reversible, hence has a unique predecessor for every directed
edge.  Thus the trajectory exits in finite time, and its terminal edge is an
outward boundary edge.

The same uniqueness of predecessors gives injectivity of $F$: if two inward
edges had the same outward terminal edge, reversing the common terminal segment
would force the two initial edges to coincide.  Since the domain and codomain
have the same cardinality, $F$ is bijective.  In the Manhattan model the local
map is a bijection from Manhattan incoming edges to Manhattan outgoing edges, so
the same argument applies to the restricted sets.
\end{proof}

\begin{lemma}\label{lem:parity}
In the Lorentz model, the number of left--right traversals is even.  Hence
$\Trav_N$ implies the existence of two left--right traversals.  In the Manhattan
model, the number of left--right Manhattan traversals equals the number of
right--left Manhattan traversals.  Hence $\Trav_N$ implies the existence of one of
each.
\end{lemma}

\begin{proof}
Consider first the Lorentz model.  Let
\[
        A=\{e\in\mathcal E_l^+:F(e)\in\mathcal E_l^-\}
\]
be the set of inward left edges whose trajectories return to the left boundary.
For $e\in A$, define $\iota(e)=\overline{F(e)}\in\mathcal E_l^+$.  By
reversibility, $\iota$ is an involution on $A$.  It has no fixed point: in the
usual half-edge representation of the Lorentz trajectories, boundary half-edges
have degree one, so one open component cannot have the same boundary half-edge as
both of its endpoints.  Hence $|A|$ is even.  Since
$|\mathcal E_l^+|=2n$ is even, the complement
$\mathcal E_l^+\setminus A$, which is precisely the set of initial edges of
left--right traversals, also has even cardinality.  If a right--left traversal
exists, its inverse is a left--right traversal, so $\Trav_N$ implies at least two
left--right traversals.

For the Manhattan model, all trajectories are counted in the Manhattan
orientation.  Let
\[
        A_M=\{e\in\mathcal E_l^+:F(e)\in\mathcal E_l^-\}.
\]
By bijectivity, $F(A_M)$ is the set of left-outward Manhattan edges used by
left--left trajectories and has cardinality $|A_M|$.  The admissible sets
$\mathcal E_l^+$ and $\mathcal E_l^-$ both have cardinality $n$.  Hence
\[
\#\{\hbox{left--right Manhattan traversals}\}
        =n-|A_M|.
\]
The right--left Manhattan traversals are exactly the trajectories, started from
$\mathcal E_r^-$, whose terminal edge lies in
$\mathcal E_l^-\setminus F(A_M)$; their number is also $n-|A_M|$.
\end{proof}

We now isolate the local modification used in the counting argument.  At a
Lorentz vertex, write $N,E,S,W$ for the four incident half-edges.  The three
possible one-site states are the pairings
\[
 (N,S)(E,W),\qquad (N,W)(S,E),\qquad (N,E)(S,W).
\]
Thus prescribing a pair of distinct half-edges determines a unique Lorentz
state.  At a Manhattan vertex, label the two incoming half-edges by $i_1,i_2$
and let $o_j$ be the straight outgoing continuation of $i_j$.  Its two possible
states are
\[
 i_1\mapsto o_1,\ i_2\mapsto o_2
 \qquad\hbox{and}\qquad
 i_1\mapsto o_2,\ i_2\mapsto o_1;
\]
toggling the site interchanges these two bijections.

Schematically, the two switches cut traversals at their contacts with $W$ and
reconnect them through an arc $P$ of the winding trajectory.  In the Lorentz
case the new path is $A_1PA_2^{-1}$; in the Manhattan case it is $A_1PB_2$.
The result is a left--left trajectory that retains both switch locations,
making the later preimage count recoverable.

\begin{lemma}[Two-site switching at a winding barrier]\label{lem:two-switch}
Fix a configuration in $D_N$.  Let $0\le a<b<N$, put
$S=T^{(n)}_{a,b}$, and let $W$ be a closed trajectory of this configuration
with $V(W)\subset S$ such that $V(W)$ separates $\{x\le a\}$ from
$\{x\ge b+1\}$.

\emph{(L) Lorentz model.}  Given two distinct left--right traversals
$L_1,L_2$ of this configuration, there are distinct contact vertices
$u_1,u_2\in V(W)$, together with marked occurrences of $u_i$ on both $L_i$
and $W$, and an arc $P$ of $W$ from $u_1$ to $u_2$ with no
internal occurrence of either contact vertex, such that changing exactly the
states at $u_1,u_2$ creates the left--left trajectory $Q=A_1PA_2^{-1}$.
Here $A_i$ is the prefix of $L_i$ ending at its first contact with $V(W)$.

\emph{(M) Manhattan model.}  Orient $W$ in the Manhattan direction.  Given a
Manhattan-oriented left--right traversal $L_1=A_1B_1$ and right--left traversal
$L_2=A_2B_2$ of this configuration, there are distinct contact vertices
$u_1,u_2$, together with marked occurrences on both the traversals and $W$, and an
oriented arc $P$ of $W$ from $u_1$ to $u_2$, with no internal occurrence of
either contact vertex, such that toggling exactly the sites $u_1,u_2$ creates a
left--left trajectory $Q=A_1PB_2$.
Here $u_1$ is the first contact of $L_1$ with $V(W)$ and $u_2$ is the last
contact of $L_2$ with $V(W)$.

In both cases the changes are genuine, both switched occurrences lie on $Q$,
$V(Q)\cap S\ne\varnothing$, and $V(Q)\subset\{x\le b\}$.
\end{lemma}

\begin{proof}
Every traversal meets $V(W)$.  Contacts are marked passages.  A traversal
cannot use a $W$-passage, since uniqueness would put it on the closed component
$W$; if $W$ uses both passages at a vertex, no other trajectory visits that
vertex.  Thus every contact below uses the unique complementary passage.

In the Lorentz case, let $u_i$ be the first contact of $L_i$ and write
$L_i=A_iB_i$.  If $u_1=u_2$, both traversals use the same complementary open
component, which has only one left--right orientation; this contradicts
$L_1\ne L_2$.  Choose consecutive $W$-occurrences labelled $u_1,u_2$,
reversing $W$ if necessary, and let $P$ be the intervening arc.  At $u_1$ use
the unique pairing joining the incoming half-edge of $A_1$ to the outgoing
half-edge of $P$; at $u_2$ join the incoming half-edge of $P$ to the outgoing
half-edge of $A_2^{-1}$.  Complementarity makes both pairings legal and
different from the old states, and gives $Q=A_1PA_2^{-1}$.

In the Manhattan case, take the first contact $u_1$ of $L_1$ and last contact
$u_2$ of $L_2$.  Equality would force both to use one complementary directed
passage, impossible for oppositely traversing Manhattan trajectories.  In the
fixed orientation of $W$, the cyclic list of $u_1$- and $u_2$-occurrences has a
transition $u_1$ to $u_2$; let $P$ be that arc.  Toggling at $u_1$ sends the
incoming half-edge of $A_1$ into $P$, and toggling at $u_2$ sends $P$ into
$B_2$.  Both toggles are genuine and give $Q=A_1PB_2$.

The first-contact prefixes, the Manhattan last-contact suffix, and the interior
of $P$ contain no other switch occurrence, so $Q$ follows the modified rule.
It meets $S$ along $P$.  No first-contact prefix can reach $\{x\ge b+1\}$, and
the Manhattan suffix cannot return there after its last contact: either event
would give a crossing of the separating strip avoiding $V(W)$.  Hence
$V(Q)\subset\{x\le b\}$.
\end{proof}

\begin{prop}[One-block no-traversal estimate]\label{prop:block}
There are constants $C_1=C_1(p)<\infty$ and $n_1=n_1(p)<\infty$ such that, if
$n\ge n_1$, $N=\lfloor C_1n^{10}\rfloor$, and
$\Prob_p(\Cross(100n,n))\ge c_0$, then
\[
        \Prob_p(\Trav_N)\le \frac14 .
\]
\end{prop}

\begin{proof}
Let
\[
        \mathcal G=\left\{\#\{1\le\ell\le L:\Wind_\ell\hbox{ occurs}\}
        \ge c_1C_1n^5\right\},
        \qquad L=\lfloor C_1n^9/1000\rfloor .
\]
By Lemma~\ref{lem:many}, $\Prob_p(\mathcal G^c)\le C_2/(C_1n)$ for all large
$n$, and, after increasing $n_1$ if necessary, every designated slice
$S_\ell$ is contained in $D_N$.

We shall define a relation $\mathcal R$ between configurations
$\omega\in\Trav_N\cap\mathcal G$ and modified configurations $\omega'$.  All
configurations are restricted to the finite set of one-site variables in $D_N$.
Whenever choices are not unique, we make them by a fixed deterministic rule, for
instance the lexicographically first edge sequence.

Fix $\omega\in\Trav_N\cap\mathcal G$ and a good slice $\ell$, and let
$W=W_\ell(\omega)$ be the chosen winding trajectory in $S_\ell$.  By
Lemma~\ref{lem:parity}, choose the required pair of traversals using the fixed
rule specified above.  Lemmas~\ref{lem:column-to-wind} and
\ref{lem:two-switch} produce a configuration $\Phi_\ell(\omega)$, differing
from $\omega$ at exactly two sites of $S_\ell$, and a left--left trajectory
$Q$ which meets $S_\ell$ and is contained in
$\{x\le 200\ell n+100n\}$.  Put
$(\omega,\Phi_\ell(\omega))\in\mathcal R$, and let $U_\ell\subset S_\ell$ be
the two switched sites.  If $\ell\ne k$, then $U_\ell\cap U_k=\varnothing$.
For $u\in U_\ell$, genuineness gives
$\Phi_\ell(\omega)(u)\ne\omega(u)=\Phi_k(\omega)(u)$, so the images are
distinct.  Each $\omega\in\Trav_N\cap\mathcal G$
has at least $c_1C_1n^5$ images under $\mathcal R$.

Let $\mu_N$ be the product measure on configurations in $D_N$.  Because
$p\in(0,1)$, every allowed one-site state has positive probability.  Since
$\omega$ and $\Phi_\ell(\omega)$ differ at exactly two sites,
\[
        \mu_N(\omega)\le C_p\,\mu_N(\Phi_\ell(\omega))
\]
for a constant $C_p=C_p(p)<\infty$.  One may take $C_p$ to be the square of the
ratio between the largest and smallest positive one-site probabilities.

We now bound the number of preimages of a fixed configuration $\omega'$.  If
$(\omega,\omega')\in\mathcal R$, then $\omega'$ contains the left--left trajectory
created by the surgery.  There are at most $2n$ choices for this trajectory,
because it is determined by its initial edge in $\mathcal E_l^+$.  Once this
trajectory is chosen, the slice used in the surgery is determined.  Indeed, if
$S_\ell=T^{(n)}_{a_\ell,b_\ell}$ is the used slice, then
$V(Q)\cap S_\ell\ne\varnothing$ and $V(Q)\subset\{x\le b_\ell\}$, while every
$S_k$ with $k>\ell$ lies in $\{x>b_\ell\}$.  Thus
$\ell=\max\{k:V(Q)\cap S_k\ne\varnothing\}$.
It remains to choose the two ordered marked passages of $Q$ at which the
switches were made.  At a fixed vertex an exiting trajectory has at most four
marked occurrences, since it cannot repeat a directed edge.  Hence there are at
most $(4|S_\ell|)^2=4^2(200n^2)^2$ ordered pairs of marked passages.
The two pre-surgery one-site
states have at most $3^2$ choices.  Once these data and $\omega'$ are fixed, at
most one predecessor configuration $\omega$ is obtained by restoring the two
states.  Therefore every
$\omega'$ has at most
\[
        2n\,(200n^2)^2\,4^2\,3^2\le 2\cdot10^7n^5
\]
preimages.

Double-counting $\mathcal R$ gives
\[
\begin{aligned}
 c_1C_1n^5 C_p^{-1}\,\Prob_p(\Trav_N\cap\mathcal G)
 &\le \sum_{(\omega,\omega')\in\mathcal R}\mu_N(\omega')  \\
 &\le \sum_{\omega'}2\cdot10^7n^5\mu_N(\omega')
      =2\cdot10^7n^5 .
\end{aligned}
\]
Thus
\[
        \Prob_p(\Trav_N)
        \le \Prob_p(\mathcal G^c)+\Prob_p(\Trav_N\cap\mathcal G)
        \le \frac{C_2}{C_1n}+\frac{2\cdot10^7C_p}{c_1C_1}.
\]
Choosing $C_1$ large enough, and then $n_1$ large enough, makes the right-hand
side at most $1/4$.
\end{proof}

\section{Block bootstrap and cylinder localization}
\label{sec:bootstrap}

Let $\Block(a,N)$ be the event that every trajectory incident to the section
$x=a+2N$ is contained in $T^{(n)}_{a,a+4N}$; it suffices to check the finitely
many directed edges incident to that section.

\begin{lemma}\label{lem:boot}
Assume that $\Prob_p(\Block(a,N))\ge1/2$ for every $a\in\Z$.  Then, for a
universal constant $C<\infty$ and all $M\ge CN$,
\[
        \Prob_p(B_n^{(M)})
        \ge 1-C\exp\left\{-\frac{M}{CN}\right\}.
\]
\end{lemma}

\begin{proof}
Following the finitely many incident directed edges until exit or repetition
shows that $\Block(a,N)$ depends only on variables in $T^{(n)}_{a,a+4N}$.
Thus each of the families $\Block(4jN,N)$ and
$\Block(-4(j+1)N,N)$, $0\le j<K:=\lfloor M/(4N)\rfloor$, is independent.
The probability that either family has no occurring block is at most
$2^{1-K}$.  One occurring block on each side prevents a trajectory from $x=0$
from reaching distance $M+1$, so
\[
        \Prob_p(B_n^{(M)})\ge1-2^{1-K},
\]
which implies the stated bound after increasing $C$.
\end{proof}

\begin{proof}[Proofs of Theorems~\ref{thm:lorentz-one-scale} and
\ref{thm:manhattan-one-scale}]
Fix either model and assume the corresponding crossing hypothesis; write $\Prob_p$ for its
law.  Let $C_1,n_1$ be as in Proposition~\ref{prop:block} and
$N=\lfloor C_1n^{10}\rfloor$.  The proposition applies to every translated
width-$N$ strip.  For Manhattan mirrors use translation $(a,0)$ when $a$ is
even and $(a,1)$ otherwise, which preserves parity.

If $\Block(a,N)$ fails, stop an escaping trajectory at its first block exit.
For a left exit, the segment after its last crossing from $x=a+N+1$ to
$x=a+N$ is a right--left traversal of $T^{(n)}_{a,a+N}$; a right exit similarly
gives a left--right traversal of $T^{(n)}_{a+3N,a+4N}$.  Therefore
\[
        \Prob_p(\Block(a,N)^c)\le\frac14+\frac14=\frac12.
\]
Lemma~\ref{lem:boot}, with $N=\lfloor C_1n^{10}\rfloor$ and constants absorbed
into $C(p)$, proves both theorems.
\end{proof}

\section{Planar confinement and unconditional consequences}
\label{sec:planar}

\begin{proof}[Proof of Proposition~\ref{prop:planar}]
Put $r=n_0$.  For $z\in(50r\Z)^2$, call $z$ bad if some trajectory starts in
$z+[-36r,36r]^2$ and exits $z+[-40r,40r]^2$.  The annulus
$[-40r,40r]^2\setminus[-36r,36r]^2$ is covered by a universal number $C_s$ of
translated rectangles of dimensions $100r\times r$.  If $z$ is bad, one of
these rectangles is crossed in its short direction.

Choose the covering rectangles in the standard four side positions of the
square annulus.  Each is either a parity-preserving translate of $R_{100r,r}$
or the image of such a translate under the diagonal reflection
$\sigma(x,y)=(y,x)$.  Since $z\in(50r\Z)^2$, the translations have even
coordinate difference and preserve the parity class of $x-y$ in the Manhattan
model.  Moreover,
\[
        \sigma_1(x,y)-\sigma_2(x,y)=y-x\equiv x-y\pmod2,
\]
so the reflection also preserves the Manhattan rule.  In the Lorentz model
these translations and the reflection preserve the law.  If a transformed
boundary orientation is opposite to the one in the definition of
$\Cross(100r,r)$, reverse the trajectory segment.  Consequently,
\[
        \Prob_p(z\text{ is bad})
        \le C_s\Prob_p(\Cross(100r,r))
        \le C_sc_0.
\]

Join two centers when their $40r$-boxes intersect.  This graph has degree at
most $8$.  Every self-avoiding path of length $j$ contains at least $j/10$
centers whose $40r$-boxes are disjoint, by a greedy selection.  The associated
badness events are independent, so a fixed self-avoiding path of length $j$ is
bad with probability at most $(C_sc_0)^{j/10}$.

If a trajectory starting in $[-1,1]^2$ reaches outside $[-m,m]^2$, record a
$50r$-grid center whose $36r$-box contains the starting point, and then record
a new center whenever the trajectory leaves the current $40r$-box.  Every
recorded center except possibly the last is bad, and consecutive centers are
adjacent.  Discard the last recorded center and loop-erase the remaining
sequence.  After increasing the universal constant $C$, this gives a
self-avoiding bad path of length at least $m/(200r)$ for $m\ge Cr$.
Therefore
\[
        \Prob_p(\Plane_m^c)
        \le\sum_{j\ge m/(200r)}8^j(C_sc_0)^{j/10}.
\]
Choose $c_0$ so that $\rho=8(C_sc_0)^{1/10}<1$.  The last sum is at most
$\rho^{\lceil m/(200r)\rceil}/(1-\rho)$, which is bounded by
$\exp\{-cm/r\}$ for $m\ge Cr$ after adjusting universal $c,C>0$.  Since
$r=n_0$, this proves the proposition.
\end{proof}

\begin{proof}[Proof of Corollary~\ref{cor:lorentz}]
Suppose that
$\Prob_p^{\mathrm L}(\Cross(100n_0,n_0))\le c_0$ for some $n_0>100$.
Proposition~\ref{prop:planar} then gives, for all sufficiently large $m$,
$\Prob_p^{\mathrm L}(\Plane_m^c)\le e^{-cm/n_0}$.
Let $E_m$ be the event that the fixed Lorentz ray from the origin considered in
\cite{kozma2015lower} reaches outside $[-m,m]^2$.  Since $\Plane_m$ confines
every trajectory starting in $[-1,1]^2$, one has $E_m\subseteq\Plane_m^c$.
The lower bound from \cite{kozma2015lower} therefore gives
\[
        \frac1{2m+1}
        \le \Prob_p^{\mathrm L}(E_m)
        \le \Prob_p^{\mathrm L}(\Plane_m^c)
        \le e^{-cm/n_0},
\]
a contradiction for sufficiently large $m$.  Thus
$\Prob_p^{\mathrm L}(\Cross(100n,n))>c_0$ for every $n>100$, and
Theorem~\ref{thm:lorentz-one-scale} proves the corollary.
\end{proof}

\begin{proof}[Proof of Corollary~\ref{cor:manhattan}]
If $\Prob_p^{\mathrm M}(\Cross(100n,n))\ge c_0$ for every sufficiently large
$n$, Theorem~\ref{thm:manhattan-one-scale} proves the result.  Otherwise fix
$n_0>100$ with $\Prob_p^{\mathrm M}(\Cross(100n_0,n_0))<c_0$ and apply
Proposition~\ref{prop:planar} at that scale.

Let $n$ now be a target half-circumference large enough that $r=n-1\ge Cn_0$.
For each vertex $(s,y)$ on a fixed section, translate its lift to $(0,0)$ when
$s-y$ is even and to $(0,1)$ otherwise; this preserves Manhattan parity.  The
inverse translate of $[-r,r]^2$ has vertical span $2r+1<2n$, so the quotient is
injective there.  For each section vertex separately, take a full planar product
environment and couple it identically to the cylinder environment on the
corresponding injectivity box; outside that box complete the planar environment
independently.  Only these separate marginal couplings are used below.

On the translated event $\Plane_r$, the projected trajectory therefore remains
in $[s-r,s+r]\times(\Z/2n\Z)$.  If $H_{s,r}$ is the event that some trajectory
from section $s$ leaves this slab, a union bound gives
\[
        \Prob_p^{\mathrm M}(H_{s,r})
        \le 2n\exp\{-c(n-1)/n_0\}.
\]
For large $n$ this is at most $1/2$.  By reindexing and reversibility it controls
all trajectories incident to the section: start from its incident directed edge
based at a section vertex and, if that edge points toward the vertex, reverse the
trajectory.  Taking $N=\lfloor n^{10}\rfloor$, the radius-$r$ slab about $s=a+2N$
lies in $T^{(n)}_{a,a+4N}$, and hence
$\Prob_p^{\mathrm M}(\Block(a,N))\ge1/2$ for every $a$.  Lemma~\ref{lem:boot}
proves the corollary.
\end{proof}

\section*{Acknowledgments}
The author thanks Lingfu Zhang and Jiaoyang Huang for several discussions.

\end{document}